\documentclass{amsart}
\usepackage{amssymb}
\usepackage{amsmath}

\newcommand{\mpn}{\medskip\par\noindent}

\newcommand{\bpn}{\bigskip\par\noindent}
\theoremstyle{definition}

\theoremstyle{remark}

\numberwithin{equation}{section}

\begin{document}
\newcommand{\Mod}[1]{\,(\text{\mbox{\rm mod}}\;#1)}
\title[On the integral of the product of several Euler polynomials]{Some properties on the integral of the product of several Euler polynomials }
\author{Taekyun Kim}
\begin{abstract}
In this paper, we study the formula for a product of two Euler polynomials. From this study, we derive some formulae for the integral of the product of two or more Euler polynomials.
\end{abstract}
\maketitle
\def\ord{\text{ord}_p}
\par\bigskip\noindent
\section{Introduction}
As is well known, the Bernoulli polynomials are given by the generating function as follows:
\begin{equation*}
\begin{split}
\frac{t}{e^t-1}e^{xt}=e^{B(x)t}=\sum_{n=0}^{\infty}B_{n}(x)\frac{t^n}{n!},\quad |t|<2\pi,
\end{split}
\end{equation*}
with the usual convention about replacing $B^n(x)$ by $B_n(x)$.
In the special case $x=0$, $B_n(0)=B_n$ are called the $n$-th  Bernoulli numbers (see [18-20]).
The constants $E_{k}$ in the Taylor series expansion
\begin{equation*}\tag{1}
\begin{split}
\frac{2}{e^t+1}=\sum_{n=0}^{\infty}E_{n}\frac{t^n}{n!},\quad |t|<\pi,
\end{split}
\end{equation*}
(cf. [2,11,12]) are known as the $k$-th Euler numbers. From the generating function of Euler numbers, we note that
\begin{equation*}\tag{2}
\begin{split}
E_{0}=1,\quad E_{n}=-\sum_{i=1}^{n}\binom{n}{i}E_{i},\quad \text{for $n \in \mathbb{N}$}.
\end{split}
\end{equation*}
The first few are $1$, $-\frac{1}{2}$, $0$, $\frac{1}{4},...,$ and $E_{2k}=0$ for $k=1,2,...$. The Euler polynomials are defined by
\begin{equation*}\tag{3}
\begin{split}
\frac{2}{e^t+1}e^{xt}=\sum_{n=0}^{\infty}E_{n}(x)\frac{t^n}{n!}=\sum_{n=0}^{\infty}\Bigg(\sum_{k=0}^{n}\binom{n}{k}E_{k}x^{n-k}\Bigg)\frac{t^n}{n!}.
\end{split}
\end{equation*}
Thus, by (3), we get
\begin{equation*}\tag{4}
\begin{split}
E_{n}(x)=\sum_{k=0}^{n}\binom{n}{k}E_{k}x^{n-k}=\big(E+x\big)^n,\quad \text{(see [2,5,11,12])},
\end{split}
\end{equation*}
with the usual convention about replacing $E^n$ by $E_n$.

From (1), (2) and (4), we have
\begin{equation*}\tag{5}
\begin{split}
E_{0}=1,\quad E_{n}(1)+E_{n}=\big(E+1\big)^n+E_{n}=0,\quad \text{if}\quad n \geq 1.
\end{split}
\end{equation*}
By the definition of Euler polynomials, we easily see that
\begin{equation*}\tag{6}
\begin{split}
\sum_{n=0}^{\infty}E_{n}(x)\frac{t^n}{n!}=\frac{2}{e^t+1}e^{xt}=\frac{2}{1+e^{-t}}e^{-t(1-x)}=\sum_{n=0}^{\infty}(-1)^{n}E_{n}(1-x)\frac{t^n}{n!}.
\end{split}
\end{equation*}
From (6), we have the reflection symmetric relation for Euler polynomials as follows:
\begin{equation*}\tag{7}
\begin{split}
 E_{n}(x)=(-1)^{n}E_{n}(1-x),\quad \text{(see [1-16])}.
\end{split}
\end{equation*}
By (4), we get
\begin{equation*}\tag{8}
\begin{split}
\frac{d}{dx} E_{n}(x)&=\sum_{l=0}^{n-1}\binom{n}{l}(n-l)x^{n-l-1}E_{l}=n\sum_{l=0}^{n-1}\binom{n-1}{l}E_{l}x^{n-1-l}
\\&=n\big(E+x\big)^{n-1}=nE_{n-1}(x).
\end{split}
\end{equation*}
Thus, by (8), we see that
\begin{equation*}\tag{9}
\begin{split}
\int_{0}^{1}E_{n}(x)dx&=\frac{1}{n+1}\int_{0}^{1}\frac{d}{dx}E_{n+1}(x)dx
\\&=\frac{1}{n+1}\big(E_{n+1}(1)-E_{n+1}\big)=-\frac{2}{n+1}E_{n+1},\quad \text{for}\quad n \in \mathbb{Z}_{+}=\mathbb{N}\cup\{0\}.
\end{split}
\end{equation*}
The gamma and beta functions are defined as the following definite integrals ($\alpha>0$, $\beta>0$):
\begin{equation*}\tag{10}
\begin{split}
\Gamma(\alpha)=\int_{0}^{\infty}e^{-t}t^{\alpha-1}dt,\quad
\text{(see [20])},
\end{split}
\end{equation*}
and
\begin{equation*}\tag{11}
\begin{split}
\mathrm{B}(\alpha,\beta)=\int_{0}^{1}t^{\alpha-1}(1-t)^{\beta-1}dt=\int_{0}^{\infty}\frac{t^{\alpha-1}}{(1+t)^{\alpha+\beta}}dt.
\end{split}
\end{equation*}
By (10) and (11), we get the following equation:
\begin{equation*}\tag{12}
\begin{split}
\Gamma(\alpha+1)=\alpha\Gamma(\alpha), \quad
\mathrm{B}(\alpha,\beta)
=\frac{\Gamma(\alpha)\Gamma(\beta)}{\Gamma(\alpha+\beta)},\quad
\text{(see [20])}.
\end{split}
\end{equation*}
In this paper we give some interesting properties of several Euler polynomials to express the integral of those polynomials from $0$
to $1$ in terms of beta and gamma functions. Finally, we derive some identities on the integral of the product of Euler polynomials.

\section{On the integral of the product of Euler polynomials}
Let us consider the integral for the product of Euler polynomials and $x^n$ as follows:
\begin{equation*}\tag{13}
\begin{split}
\int_{0}^{1}y^nE_{n}(x+y)dy=\sum_{l=0}^{n}\binom{n}{l}E_{n-l}(x)\int_{0}^{1}y^{n+l}dy
=\sum_{l=0}^{n}\frac{\binom{n}{l}E_{n-l}(x)}{n+l+1}.
\end{split}
\end{equation*}
On the other hand, by (7), we get
\begin{equation*}\tag{14}
\begin{split}
\int_{0}^{1}y^nE_{n}(x+y)dy&=(-1)^n\int_{0}^{1}y^{n}E_{n}(1-(x+y))dy
\\&=(-1)^n\sum_{l=0}^{n}\binom{n}{l}E_{n-l}(-x)\int_{0}^{1}y^{n}(1-y)^ldy
\\&=\sum_{l=0}^{n}\binom{n}{l}(-1)^lE_{n-l}(1+x)\mathrm{B}(n+1,l+1)
\\&=\sum_{l=0}^{n}\binom{n}{l}(-1)^lE_{n-l}(1+x)\frac{\Gamma(n+1)\Gamma(l+1)}{\Gamma(n+l+2)}.
\end{split}
\end{equation*}
Therefore, by (13) and (14), we obtain the following theorem.
\\

{\bf Theorem 1.}  {\it For $n \in\mathbb{Z}_{+}=\mathbb{N}\cup\{0\}$, we have}
\begin{equation*}
\begin{split}
\sum_{l=0}^{n}\frac{\binom{n}{l}E_{n-l}(x)}{n+l+1}=\sum_{l=0}^{n}(-1)^l\Bigg(\frac{E_{n-l}(1+x)}{n+l+1}  \Bigg)\frac{\binom{n}{l}}{\binom{n+l}{l}}.
\end{split}
\end{equation*}
 {\it In particular, $x=0$,}
\begin{equation*}
\begin{split}
\sum_{l=0}^{n}\frac{\binom{n}{l}E_{n-l}}{n+l+1}=(-1)^n\sum_{l=0}^{n}\frac{E_{n-l}}{n+l+1} \frac{\binom{n}{l}}{\binom{n+l}{l}}.
\end{split}
\end{equation*}
\\

Let $n \in \mathbb{N}$ with $n \geq 3$. Then, by (9), we see that
\begin{equation*}\tag{15}
\begin{split}
&\int_{0}^{1}y^nE_{n}(x+y)dy=\frac{E_{n}(x+1)}{n+1} -\frac{n}{n+1}\int_{0}^{1}y^{n+1}E_{n-1}(x+y)dy
\\&=\frac{E_{n}(x+1)}{n+1}-\frac{E_{n-1}(x+1)}{n+1}\frac{n}{n+2}+(-1)^{2}\frac{n(n-1)}{(n+1)(n+2)}\int_{0}^{1}y^{n+2}E_{n-2}(x+y)dy
\\&=\frac{E_{n}(x+1)}{n+1}-\frac{nE_{n-1}(x+1)}{(n+1)(n+2)}+(-1)^{2}\frac{n(n-1)E_{n-2}(x+1)}{(n+1)(n+2)(n+3)}\\& \quad\quad\quad +(-1)^{3}\frac{n(n-1)(n-2)}{(n+1)(n+2)(n+3)}\int_{0}^{1}y^{n+3}E_{n-3}(x+y)dy.
\end{split}
\end{equation*}
Continuing this process, we obtain the following equation:
\begin{equation*}\tag{16}
\begin{split}
&\int_{0}^{1}y^nE_{n}(x+y)dy
\\&=\frac{E_{n}(x+1)}{n+1} +\sum_{l=2}^{n-1}\frac{n(n-1)\cdots(n-l+2)(-1)^{l-1}}{(n+1)(n+2) \cdots (n+l)}E_{n-l+1}(1+x)
\\& \quad\quad\quad+(-1)^{n-1}\frac{n(n-1)(n-2)\cdots 2}{n(n+1)\cdots(2n-1)}\int_{0}^{1}y^{2n-1}E_{1}(x+y)dy
\\&=\frac{E_{n}(x+1)}{n+1} -\sum_{l=2}^{n-1}\frac{n(n-1)\cdots(n-l+2)(-1)^{l-1}}{(n+1)(n+2) \cdots (n+l)}E_{n-l+1}(1+x)
\\& \quad\quad\quad+(-1)^{n-1}\frac{n!}{(n+1)(n+2)\cdots2n}\Big(E_{1}(x+1) -\frac{1}{2n+1} \Big).
\end{split}
\end{equation*}
Therefore, by (14) and (16), we obtain the following theorem.
\\

{\bf Theorem 2.}  {\it For $n \in \mathbb{N}$ with $n \geq 3$, we have}
\begin{equation*}
\begin{split}
\sum_{l=0}^{n}(-1)^l\frac{E_{n-l}(x+1)}{n+l+1} \frac{\binom{n}{l}}{\binom{n+l}{l}}
&=\frac{E_{n}(x+1)}{n+1} +\sum_{l=2}^{n-1}\frac{\binom{n}{l}(n-l+2)(-1)^{l-1}}{\binom{n+l}{l}}E_{n-l+1}(1+x)
\\&\quad\quad+\frac{(-1)^{n-1}}{\binom{2n}{n}}\Big(E_{1}(x+1) -\frac{1}{2n+1} \Big).
\end{split}
\end{equation*}
 {\it In the special case, $x=0$,}
\begin{equation*}
\begin{split}
\sum_{l=0}^{n}\frac{E_{n-l}}{n+l+1} \frac{\binom{n}{l}}{\binom{n+l}{l}}&=\frac{E_{n}}{n+1} +\sum_{l=2}^{n-1}\frac{\binom{n}{l}(n-l+2)}{\binom{n+l}{l}}E_{n-l+1}
+\frac{1}{\binom{2n}{n}}\Big(\frac{1}{2}+\frac{1}{2n+1} \Big).
\end{split}
\end{equation*}
\\

For $n \in \mathbb{N}$, we have
\begin{equation*}\tag{17}
\begin{split}
\int_{0}^{1}y^nE_{n}(x+y)dy
=\frac{E_{n}(x)}{n+1}-\frac{n}{n+1}\int_{0}^{1}y^{n-1}E_{n+1}(x+y)dy.
\end{split}
\end{equation*}
By (17), we get
\begin{equation*}\tag{18}
\begin{split}
\int_{0}^{1}y^nE_{n}(x+y)dy
=\frac{E_{n}(x)}{n+1}-\frac{n}{n+1}\sum_{l=0}^{n+1}\binom{n+1}{l}E_{n+1-l}(-x)(-1)^{n+1}\mathrm{B}(n,l+1).
\end{split}
\end{equation*}
Therefore, by (14) and (18), we obtain the following theorem.
\\

{\bf Theorem 3.}  {\it For $n \in \mathbb{N}$,  we have}
\begin{equation*}
\begin{split}
\frac{E_{n+1}(x)}{n+1}=\frac{1}{n+1}\sum_{l=0}^{n+1}\frac{\binom{n+1}{l}}{\binom{n+l}{l}}E_{n+1-l}(1+x)(-1)^{l}-\sum_{l=0}^{n}\frac{E_{n-l}(1+x)}{n+l+1}(-1)^{l}\frac{\binom{n}{l}}{\binom{n+l}{l}}
\end{split}
\end{equation*}
 {\it In particular, $x=0$,}
\begin{equation*}
\begin{split}
-\frac{E_{n+1}}{n+1}=\frac{1}{n+1}\sum_{l=0}^{n+1}\frac{\binom{n+1}{l}}{\binom{n+l}{l}}E_{n+1-l}+\sum_{l=0}^{n}\frac{\binom{n}{l}}{\binom{n+l}{l}}\frac{E_{n-l}}{n+l+1}.
\end{split}
\end{equation*}
\\

A definite integral for the multiplication of two Euler polynomials can be given by the following relation:
\begin{equation*}\tag{19}
\begin{split}
\int_{0}^{1}E_{n}(x)E_{m}(x)dx
&=\sum_{l=0}^{n}\binom{n}{l}E_{l}(-1)^{m}\sum_{k=0}^{m}\binom{m}{k}E_{k}  \int_{0}^{1}x^{n-l}(1-x)^{m-k}dx
\\&=\sum_{l=0}^{n}\sum_{k=0}^{m}\binom{n}{l}\binom{m}{k}(-1)^{m}E_{l}E_{k}\mathrm{B}(n-l+1,m-k+1).
\end{split}
\end{equation*}
Let $m,n \in \mathbb{N}_{+}$ with $m \geq 1$. Then we see that
\begin{equation*}\tag{20}
\begin{split}
\int_{0}^{1}E_{n}(x)E_{m}(x)dx
&=-\frac{m}{n+1}\int_{0}^{1}E_{n+1}(x)E_{m-1}(x)dx
\\&=(-1)^{2}\frac{m(m-1)}{(n+1)(n+2)}\int_{0}^{1}E_{n+2}(x)E_{m-2}(x)dx.
\end{split}
\end{equation*}
Continuing this process, we have
\begin{equation*}\tag{21}
\begin{split}
\int_{0}^{1}E_{n}(x)E_{m}(x)dx
=(-1)^{m-1}\frac{m(m-1)\cdots 2}{(n+1)(n+2)\cdots (n+m-1)}\int_{0}^{1}E_{n+m-1}(x)E_{1}(x)dx.
\end{split}
\end{equation*}
It is easy to show that
\begin{equation*}\tag{22}
\begin{split}
\int_{0}^{1}E_{n+m-1}(x)E_{1}(x)dx=-\frac{1}{m+n}\int_{0}^{1}E_{m+n}(x)=\frac{2E_{n+m-1}}{(m+n)(m+n+1)}.
\end{split}
\end{equation*}
Thus, by (20), (21) and (22), we get
\begin{equation*}\tag{23}
\begin{split}
\int_{0}^{1}E_{n}(x)E_{m}(x)dx
&=(-1)^{m+1}\frac{2m!E_{n+m+1}}{(n+1)(n+2)\cdots (n+m)(n+m+1)}\\&=2(-1)^{m+1}\frac{m!n!}{(n+m)!}\frac{E_{n+m+1}}{n+m+1}.
\end{split}
\end{equation*}
Therefore, by (19) and (23), we obtain the following theorem.
\\

{\bf Theorem 4.}  {\it For $m,n \in \mathbb{Z}_{+}$ with $m \geq 1$,  we have}
\begin{equation*}
\begin{split}
-2\frac{E_{n+m+1}}{\binom{n+m}{n}(n+m+1)}=\sum_{l=0}^{n}\sum_{k=0}^{m}\binom{n}{l}\binom{m}{k}E_{l}E_{k}\frac{\Gamma(n-l+1)\Gamma(m-k+1)}{\Gamma(n+m-l-k+2)}.
\end{split}
\end{equation*}
 {\it Moveover}
\begin{equation*}
\begin{split}
E_{n+m+1}=-\frac{1}{2}\sum_{l=0}^{n}\sum_{k=0}^{m}\frac{\binom{n}{l}\binom{m}{k}\binom{n+m}{n}}{\binom{n+m-l-k}{n-l}}\Bigg(\frac{(n+m+1)E_{l}E_{k}}{n+m-l-k+1}\Bigg).
\end{split}
\end{equation*}
\\
From (3), we note that
\begin{equation*}\tag{24}
\begin{split}
\sum_{m,n=0}^{\infty}& \Big(mE_{m-1}(x)E_{n}(x)+nE_{m}(x)E_{n-1}(x)  \Big)\frac{t^m}{m!}\frac{s^n}{n!}
\\&=\frac{d}{dx}\Bigg(\frac{4e^{(s+t)x}}{(e^t+1)(e^s+1)}\Bigg)=(s+t)\frac{4}{(e^t+1)(e^s+1)}e^{(s+t)x}
\\&=\Bigg(\frac{s+t}{e^{(s+t)-1}}e^{(s+t)x}\Bigg)\Bigg(4-\frac{4}{e^t+1}-\frac{4}{e^s+1}  \Bigg)
\\&=\Bigg(\sum_{l=0}^{\infty}B_{l}(x)\frac{(s+t)^l}{l!}\Bigg)\Bigg(4-2\sum_{r=0}^{\infty}E_{r}\frac{t^r}{r!}-2\sum_{r=0}^{\infty}E_{r}\frac{s^r}{r!}  \Bigg)
\\&=\Bigg(\sum_{m,n=0}^{\infty}B_{m+n}(x)\frac{t^ms^n}{m!n!}\Bigg)\Bigg(-2\sum_{r=0}^{\infty}\frac{E_{2r+1}}{(2r+1)!}\Big( t^{2r+1}+s^{2r+1}\Big) \Bigg)
\\&=\sum_{m,n=0}^{\infty}\Bigg( -2\sum_{r=0}^{\infty}\frac{E_{2r+1}}{(2r+1)!} \frac{B_{m+n}(x)}{m!n!} \Big( t^{2r+m+1}s^n+s^{2r+1+n}t^m\Big)        \Bigg)
\\&=\sum_{m,n=0}^{\infty}\Bigg( -2\sum_{r=0}^{\infty}\frac{E_{2r+1}}{(2r+1)!}  \Big( B_{m-2r-1+n}(x)\frac{t^ms^n m!}{(m-2r-1)!n!m!}  \Big)
\\&\quad\quad\quad\quad\quad\quad+\frac{ B_{m+n-2r-1}(x)}{(n-2r-1)!m!}s^n t^m\frac{n!}{n!}    \Bigg)
\\&=\sum_{m,n=0}^{\infty}\Bigg( -2\sum_{r=0}^{\infty}E_{2r+1}B_{m-2r-1+n}(x)\Bigg( \binom{m}{2r+1}+\binom{n}{2r+1} \Bigg)\Bigg) \frac{t^ms^n }{m!n!}.
\end{split}
\end{equation*}
By comparing coefficients on the both sides in (24), we obtain the following theorem.
\\

{\bf Theorem 5.}  {\it For $m,n \in \mathbb{N}$,  we have}
\begin{equation*}
\begin{split}
mE_{m-1}(x)E_{n}(x)+nE_{m}(x)E_{n-1}(x)=-2\sum_{r=0}^{\infty}E_{2r+1}B_{m-2r-1+n}(x)\Bigg( \binom{m}{2r+1}+\binom{n}{2r+1} \Bigg).
\end{split}
\end{equation*}
 \\
From (8) and (9), we note that
\begin{equation*}\tag{25}
\begin{split}
\frac{d}{dx}\Big( E_{m}(x)E_{n}(x)\Big)&=mE_{m-1}(x)E_{n}(x)+nE_{m}(x)E_{n-1}(x)
\\&=-2\sum_{r=0}^{\infty}E_{2r+1}B_{m-2r-1+n}(x)\Bigg( \binom{m}{2r+1}+\binom{n}{2r+1} \Bigg).
\end{split}
\end{equation*}
Thus, by (25), we get
\begin{equation*}\tag{26}
\begin{split}
 E_{m}(x)E_{n}(x)&=\int_{}^{} \frac{d}{dx}\Big( E_{m}(x)E_{n}(x)\Big)dx
\\&=-2\sum_{r=0}^{\infty}E_{2r+1}\frac{B_{m+n-2r}(x)}{m+n-2r}\Bigg( \binom{m}{2r+1}+\binom{n}{2r+1} \Bigg)+C,
\end{split}
\end{equation*}
where $C$ is some constant.

For $m+n\geq 2$, we have
\begin{equation*}\tag{27}
\begin{split}
\int_{0}^{1}E_{m}(x)E_{n}(x)dx=C.
\end{split}
\end{equation*}
By (23) and (27), we get
\begin{equation*}\tag{28}
\begin{split}
C=2(-1)^{m+1}\frac{m!n!}{(n+m)!}\frac{E_{n+m+1}}{n+m+1}=(-1)^{m+1}\frac{2}{\binom{n+m}{n}}\frac{E_{n+m+1}}{n+m+1}.
\end{split}
\end{equation*}
Therefore, by (26) and (28), we obtain the following theorem.
\\

{\bf Theorem 6.}  {\it For $m,n \in \mathbb{N}$ with $m+n \geq 2$,  we have}
\begin{equation*}
\begin{split}
&E_{m}(x)E_{n}(x)\\&=-2\sum_{r=0}^{\infty}\Bigg( \binom{m}{2r+1}+\binom{n}{2r+1} \Bigg)\frac{E_{2r+1}B_{m+n-2r}(x)}{m+n-2r}+2(-1)^{m+1}\frac{E_{n+m+1}}{\binom{n+m}{n}(n+m+1)}.
\end{split}
\end{equation*}
\\
Note that
\begin{equation*}
\begin{split}
&\int_{0}^{1}B_{m+n-2r}(x)E_{p}(x)dx
=\Bigg[ \frac{B_{m+n-2r+1}(x)E_{p}(x)}{m+n-2r+1}           \Bigg]_{0}^{1}-p\int_{0}^{1}\frac{B_{m+n-2r+1}(x)E_{p-1}(x)}{m+n-2r+1}dx
\\&=-2\frac{B_{m+n-2r+1}E_{p}}{m+n-2r+1} -p\Bigg[ \frac{B_{m+n-2r+2}(x)E_{p-1}(x)}{(m+n-2r+1)(m+n-2r+2)} \Bigg]_{0}^{1}
\\& \quad\quad\quad\quad\quad\quad\quad\quad\quad+p(p-1)(-1)^2 \int_{0}^{1} \frac{B_{m+n-2r+2}(x)E_{p-2}(x)}{(m+n-2r+1)(m+n-2r+2)}dx.
\end{split}
\end{equation*}
Continuing this process, we obtain
\begin{equation*}\tag{29}
\begin{split}
&\int_{0}^{1}B_{m+n-2r}(x)E_{p}(x)dx
\\&=2p! \sum_{l=0}^{p-1}\frac{B_{m+n-2r+l}E_{p-l+1}(-1)^l}{(m+n-2r+1)\cdots(m+n-2r+l)}+\frac{p!(-1)^{p-1}\int_{0}^{1}B_{m+n-2r+p-1}(x)E_{1}(x)dx}{(m+n-2r+1)\cdots(m+n-2r+p-1)}
\\&=2p! \sum_{l=0}^{p-1}\frac{B_{m+n-2r+l}E_{p-l+1}(-1)^l}{(m+n-2r+1)\cdots(m+n-2r+l)}+\frac{2p!(-1)^{p}B_{m+n-2r+p}E_{1}}{(m+n-2r+1)\cdots(m+n-2r+p)}
\\&=2p! \sum_{l=0}^{p}\frac{B_{m+n-2r+l}E_{p-l+1}(-1)^l}{(m+n-2r+1)\cdots(m+n-2r+l)}
=2p! \sum_{l=0}^{p}\frac{B_{m+n-2r+l}E_{p-l+1}}{\binom{n+m-2r+l}{l}}\frac{(-1)^l}{l!}.
\end{split}
\end{equation*}
From Theorem 6, we can derive the following equation (30). For $m,n,p \in \mathbb{N}$,  we have
\begin{equation*}\tag{30}
\begin{split}
&\int_{0}^{1} E_{m}(x)E_{n}(x)E_{p}(x)dx=\int_{0}^{1} E_{p}(x)(E_{m}(x)E_{n}(x))dx
\\&=-2\sum_{r=0}^{\infty}\Bigg( \binom{m}{2r+1}+\binom{n}{2r+1} \Bigg)\frac{E_{2r+1}}{m+n-2r} \int_{0}^{1}B_{m+n-2r}(x)E_{p}(x)dx
\\&\quad\quad\quad+2(-1)^{m+1}\frac{E_{n+m+1}}{\binom{n+m}{n}(n+m+1)}\int_{0}^{1}E_{p}(x)dx.
\end{split}
\end{equation*}
By (29) and (30), we get
\begin{equation*}\tag{31}
\begin{split}
&\int_{0}^{1} E_{m}(x)E_{n}(x)E_{p}(x)dx
\\&=-4p!\sum_{r=0}^{\infty}\Bigg( \binom{m}{2r+1}+\binom{n}{2r+1} \Bigg)\frac{E_{2r+1}}{n+m-2r}\sum_{l=1}^{p}\frac{B_{m+n-2r+l}E_{p-l+1}}{\binom{n+m-2r+l}{l}}\frac{(-1)^l}{l!}
\\&\quad\quad\quad+4\frac{(-1)^{m}E_{n+m+1}E_{p+1}}{\binom{n+m}{n}(n+m+1)(p+1)}.
\end{split}
\end{equation*}
On the other hand,
\begin{equation*}\tag{32}
\begin{split}
&\int_{0}^{1} E_{m}(x)E_{n}(x)E_{p}(x)dx
\\&=(-1)^{p}\sum_{l=0}^{m}\sum_{j=0}^{n}\sum_{k=0}^{p}\binom{m}{l}\binom{n}{j}\binom{p}{k}E_{m-l}E_{n-j}E_{k}\int_{0}^{1}x^{l+j}(1-x)^{p-k}dx
\\&=(-1)^{p}\sum_{l=0}^{m}\sum_{j=0}^{n}\sum_{k=0}^{p}\binom{m}{l}\binom{n}{j}\binom{p}{k}E_{m-l}E_{n-j}E_{k}\mathrm{B}(l+j+1,p-k+1)
\\&=(-1)^{p}\sum_{l=0}^{m}\sum_{j=0}^{n}\sum_{k=0}^{p}\binom{m}{l}\binom{n}{j}\binom{p}{k}E_{m-l}E_{n-j}E_{k}\frac{(l+j)!(p-k)!}{(l+j+p-k+1)!}.
\end{split}
\end{equation*}
Therefore, by (31) and (32), we obtain the following equation:
\begin{equation*}
\begin{split}
&-4p!\sum_{r=0}^{\infty}\Bigg( \binom{m}{2r+1}+\binom{n}{2r+1} \Bigg)\frac{E_{2r+1}}{n+m-2r}\sum_{l=1}^{p}\frac{B_{m+n-2r+l}E_{p-l+1}}{\binom{n+m-2r+l}{l}}\frac{(-1)^l}{l!}
\\&+4\frac{(-1)^{m}E_{n+m+1}E_{p+1}}{\binom{n+m}{n}(n+m+1)(p+1)}
=(-1)^{p}\sum_{l=0}^{m}\sum_{j=0}^{n}\sum_{k=0}^{p}\frac{\binom{m}{l}\binom{n}{j}\binom{p}{k}E_{m-l}E_{n-j}E_{k}}{\binom{l+j+p-k}{l+j}(l+j+p-k+1)},
\end{split}
\end{equation*}
where $m,n,p \in \mathbb{N}$.

\par\bigskip
ACKNOWLEDGEMENTS. The author would like to express his sincere gratitude to referees for their valuable comments and suggestions.
%
%

\vspace{30mm}

\par\bigskip
\begin{center}\begin{large}

{\sc References}
\end{large}\end{center}
\par
\begin{enumerate}

\item[{[1]}] M. Abramowitz, I. A. Stegun, {\it Handbook of mathematical functions with formulas, graphs, and mathematical tables},
National Bureau of Standards Applied Mathematics Series, 55
(1964).

\item[{[2]}] S. Araci, D. Erdal, J. Seo, {\it A study on the fermionic $p$-adic $q$-integral on $\mathbb{Z}_{p}$ associated with weighted  $q$-Bernstein  $q$-Genocchi polynomials},
Abstract and Applied Analysis, 2011(2011), Article ID {\bf 649248}, 10 pages.

\item[{[3]}] A. Bayad, J. Chikhi, {\it Non linear recurrences for Apostol-Bernoulli-Euler number of higher order},
Adv. Stud. Contemp. Math. \textbf{22}(2012), no.1, 1--6.

\item[{[4]}] A. Bayad, T. Kim, B. Lee, S.-H. Rim, {\it Some identities on the  Bernstein polynomials associated with $q$-Euler  polynomials},
Abstract and Applied Analysis, 2011(2011), Article ID {\bf 294715}, 10 pages.

\item[{[5]}]  L. Carlitz, {\it  Note on the integral of the product of several Bernoulli polynomials},  J. Londen Math. Soc. 34 (1959), 361--363.

\item[{[6]}]  M. Cenkci,  Y. Simsek,  V. Kurt, {\it Multiple two-variable $q$-$L$-function and its behavior at $s=0$},
 Russ. J. Math. Phys. 15 (2008), 447--459.

\item[{[7]}]  J. Choi, D. S. Kim, T. Kim, Y.-H. Kim,
 {\it A note on some identities of Frobenius-Euler numbers and polynomials},
 Int. J. Math. Math. Sci. 2012(2012), Article ID 861797, 7 pages.

 \item[{[8]}]  T. Kim, {\it  Some identities on the $q$-Euler polynomials of
higher order and $q$-Stirling numbers by the fermionic $p$-adic
integral on $\Bbb Z_p$}, Russ. J. Math. Phys. 16 (2009),
484--491.

 \item[{[9]}]  T. Kim, {\it Symmetry of power sum polynomials and multivariate fermionic $p$-adic invariant
integral on $\Bbb Z_p$}, Russ. J. Math. Phys. 16 (2009),
93--96.

\item[{[10]}] T. Kim,  {\it A note on $q$-Bernstein polynomials}, Russ. J. Math. Phys. 18 (2011), 73--82.

\item[{[11]}] T. Kim,  {\it A note on the Euler numbers and polynomials}, Adv. Stud. Contemp. Math. \textbf{17}(2008), 131--136.

\item[{[12]}] T. Kim,  {\it Some identities on the $q$-integral representation of the product
of several $q$-Bernstein-type polynomials}, Abstract and Applied
Analysis, \textbf{2011}(2011), Article ID 634675, 11 pages.

\item[{[13]}] H. Ozden, I. N. Cangul, Y. Simsek, {\it Remarks on $q$-Bernoulli numbers associated with Daehee numbers},
Adv. Stud. Contemp. Math. \textbf{18} (2009), 41--48.

\item[{[14]}] H. Ozden, I. N. Cangul, Y. Simsek, {\it Multivariate interpolation functions of higher-order $q$-Euler numbers and their applications},
Abstr.  Appl. Anal. 2008(2008), Art. ID. {\bf 390857}, 16 pages.

\item[{[15]}] C. S. Ryoo, {\it Some identities of the twisted $q$-Euler numbers and polynomials associated with $q$-Bernstein  polynomials}, Proc. Jangjeon
Math. Soc. \textbf{14} (2011), 239--348.

\item[{[16]}] C. S. Ryoo, {\it Some relations between twisted $q$-Euler numbers and Bernstein  polynomials},
Adv. Stud. Contemp. Math. \textbf{21}(2011), 217--223.

\item[{[17]}] Y. Simsek, {\it Generating functions of the twisted Bernoulli numbers and polynomials associated with their interpolation functions},
Adv. Stud. Contemp. Math. \textbf{16} (2008), 251--278.

\item[{[18]}] Y. Simsek, {\it Theorems on twisted $L$-function and twisted Bernoulli numbers},
Adv. Stud. Contemp. Math. \textbf{11} (2005), 205--218.

\item[{[19]}]  Y. Simsek,  {\it Special functions related to Dedekind-type DC-sums and their applications},
 Russ. J. Math. Phys. 17 (2010), 495--508.

\item[{[20]}]  D. G. Zill, M. R. Cullen, {\it Advanced Engineering Mathematics},
 Jones and Bartlett Publishers, Inc., 2006.

\end{enumerate}

\vspace{5mm}

\mpn { \bpn {\small Taekyun {\sc Kim} \mpn Department of
Mathematics,  Kwangwoon University, Seoul 139-701, Republic of
Korea,  {\it E-mail:}\ {\sf tkkim@kw.ac.kr} }

\end{document}